\documentclass[preprint,12pt]{elsarticle}
\journal{International Journal of Electrical Power \& Energy Systems}
\usepackage[top=3cm,bottom=3cm,left=3cm,right=3cm]{geometry}
\usepackage[cmex10]{amsmath}
\usepackage{comment}
\usepackage{ifpdf} 
\usepackage{graphicx} 
\usepackage{epstopdf}
\usepackage{multirow}
\usepackage[font=small,labelfont=rm]{caption}
\usepackage{tabularx,ragged2e}
\usepackage{algorithm}
\usepackage{hyperref}
\usepackage{algorithmicx}
\usepackage{siunitx}
\usepackage{algpseudocode}
\usepackage{graphics}
\usepackage{etoolbox}
\usepackage{nomencl}
\makenomenclature
\usepackage{amssymb}
\usepackage{lineno}
\usepackage{bm}
\usepackage[dvipsnames]{xcolor}

\newcommand{\norm}[1]{\left\lVert#1\right\rVert} % For the SOC

\usepackage{tikz}

\usepackage{pgfplots}
\pgfplotsset{width=5cm,compat=1.9}
\usepgfplotslibrary{external}
\usepgfplotslibrary{statistics}

\usepackage{etoolbox}
\makeatletter
\patchcmd{\ps@pprintTitle}% <cmd>
  {submitted to}% <search>
  {accepted for publication in the}% <replace>
  {}{}% <succes><failure>
\makeatother

\begin{document}

\begin{frontmatter}

\title{Chance-constrained allocation of UFLS candidate feeders under high penetration of distributed generation}
\author[1]{Luis Badesa\corref{cor1}\fnref{fn1}}
\ead{luis.badesa@upm.es}
\author[2]{Cormac O'Malley\fnref{fn1}}
\ead{c.omalley19@imperial.ac.uk}
\author[3]{Mar{\'i}a Parajeles}
\ead{parajeles_herrera@eeh.ee.ethz.ch}
\author[2]{Goran Strbac}
\ead{g.strbac@imperial.ac.uk}

\cortext[cor1]{Corresponding author}
\address[1]{Technical University of Madrid (UPM), Ronda de Valencia 3, 28012 Madrid,
Spain}
\address[2]{Imperial College London, South Kensington, London SW7 2BU, United Kingdom}
\address[3]{ETH Zurich, R{\"a}mistrasse 101, Z{\"u}rich 8092, Switzerland}

\fntext[fn1]{Luis Badesa and Cormac O'Malley contributed equally to this work.}

\begin{abstract}
Under-Frequency Load Shedding (UFLS) schemes are the last resort to contain a frequency drop in the grid by disconnecting part of the demand. The allocation methods for selecting feeders that would contribute to the UFLS scheme have traditionally relied on the fact that electric demand followed fairly regular patterns, and could be forecast with high accuracy. However, recent integration of Distributed Generation (DG) increases the uncertainty in net consumption of feeders which, in turn, requires a reformulation of UFLS-allocation methods to account for this uncertainty. In this paper, a chance-constrained methodology for selecting feeders is proposed, with mathematical guarantees for the disconnection of the required amount of load with a certain pre-defined probability. The correlation in net-load forecasts among feeders is explicitly considered, given that uncertainty in DG power output is driven by meteorological conditions with high correlation across the network. Furthermore, this method is applicable either to systems with conventional UFLS schemes (where relays measure local frequency and trip if this magnitude falls below a certain threshold), or adaptive UFLS schemes (where relays are triggered by control signals sent in the few instants following a contingency). Relevant case studies demonstrate the applicability of the proposed method, and the need for explicit consideration of uncertainty in the UFLS-allocation process.
\end{abstract}

 \begin{highlights}
     \item A chance-constrained methodology for selecting feeders contributing towards UFLS is proposed. 
     \item Uncertainty and correlation driven by Distributed Generation are explicitly considered.
     \item A distributionally-robust option is given, requiring minimal information on uncertainty.
     \item The method applies to several communication requirements for the operation of the power system.
 \end{highlights}


\end{frontmatter}

\section{Introduction} \label{sec:Intro}

The deployment of renewable energy sources (RES) in electricity grids across the world has experienced a very significant increase in recent years, not only on the transmission level, but also on the distribution side: more than $30\%$ of the solar photovoltaic capacity installed is connected to the distribution system \cite{IEARenewables20192024}. The rise in this distributed generation (DG) brings several challenges to system operators, given the limited observability and controllability they have over DG. Furthermore, given that most RES are non-synchronous generators, the level of inertia in modern grids is also significantly reduced, leading to more severe frequency variations \cite{FrequencyVolatility}. Under-Frequency Load Shedding (UFLS) then becomes increasingly relevant, as the last-resource containment mechanism for keeping the system in balance in the event of a large power infeed loss.

A UFLS scheme disconnects feeders from the system to restore demand-generation balance and thus arrest frequency decline. However, this has the negative consequence of interrupting supply to some consumers. In traditional UFLS schemes, protection relays measure their local frequency and trip when this magnitude breaches a predefined trigger level \cite{Rowland}. The combination of relays to be disconnected at each trigger level is chosen to disconnect a set percentage of national demand. This technique is highly reliable and widely deployed, but has the disadvantage of indiscriminately shedding the same amount of load for a wide range of outage sizes and inertia levels, thus leading to costly over-shedding of demand.

The improved design of UFLS schemes to contain frequency with the least load shed possible has been an area of intense research \cite{Haesalhelou2020,Sanaye-Pasand2009,Terzija2006,Rudez2016}. Schemes utilise a post-outage frequency nadir prediction method to calculate online the optimal amount of load to be shed \cite{Haesalhelou2020}. Works such as \cite{Sanaye-Pasand2009} assume the outage size to be known and use the single-machine equivalent swing equation \cite{KundurBook} to predict the nadir. When the outage size is not known, it can be estimated from initial rate-of-change-of-frequency measurements \cite{Terzija2006}, while approaches such as \cite{Rudez2016} estimate the nadir directly from frequency measurements. 

These methods all focus on determining the amount of load to be shed. However, because DG make the net-power flow through feeders stochastic, it will no longer be a simple task to map a desired amount of demand to be shed into a given choice of feeders. Shedding less load than desired will hinder the efficacy of the proposed UFLS schemes and jeopardise system security. Indeed, in extreme cases, DG can make some feeders power exporters at times, in which case their disconnection will exacerbate the problem. Solving this second stage relay-allocation problem is vital to the success of any UFLS scheme, and currently no robust method to deal with the uncertainty exists.

A few limited heuristic methods represent the current state of the art for the feeder allocation problem \cite{Li2017,Atighechi2018, Chandra2021,DeBoeck2018,Das2017,IET_UFLS_clustering}. References \cite{Li2017}, \cite{Atighechi2018} and \cite{Chandra2021} analyse the effects of disconnecting different feeders on the system's dynamic performance and voltage stability. The prioritisation of feeders to be disconnected is based on certain voltage and frequency stability indexes, linked to active power changes. The online computation of indexes requires real-time measurements, which make these methods highly dependent on measurement and communication devices to perform the UFLS allocation among feeders. The authors of \cite{DeBoeck2018} and \cite{Das2017} explore the prioritisation of feeders through the concept of generation-to-consumption ratio, i.e.,~the total amount of load disconnected compared to the net-power flow disconnected. This metric is used to reduce disruption to customers, as feeders are  heuristically classified based on the values of this ratio, prioritising the disconnection of the most load-dense circuits. This classification leads to some of the consumers never being considered for the UFLS scheme, and the selection of feeders is still dependent on real-time power flow measurements and DG output estimation, therefore no uncertainty is modelled. Finally, the authors in \cite{IET_UFLS_clustering} compare current UFLS-allocation methods used in France, based on simple historical measurements of consumption through each feeder, with their proposed method based on clustering feeders that could jointly contribute to minimise deviation of actual load shed from the desired volume of load shed. The clustering is computed using also a time-series of historical data.

None of the previous work, to the best of our knowledge, has addressed the allocation of load shedding while accounting for uncertainty, in a formulation that mathematically guarantees the disconnection of the required load. The core contribution of this paper is to propose a method that translates a desired load to be shed into a set of relays that guarantee delivery under DG uncertainty, with a user-specified level of risk for under-delivery in extreme scenarios. 

A strength of our proposed method is that it is useful in combination with all types of UFLS schemes, although its application varies depending on when the decision on which relays will be disconnected is made (offline before the outage \cite{Rowland} or online during the transient period following an outage \cite{Sanaye-Pasand2009,Terzija2006,Rudez2016}). For the first case \cite{Rowland}, where UFLS relays are triggered in a decentralised manner through local frequency measurements, the net-load through the feeders will always be uncertain at the time of outage. This is because the exact DG generation can not be exactly forecast ahead of time and thus the exact net-load through feeders cannot be known. However, accurate forecasts can produce information on the expected DG output value, which our method utilises alongside the forecast error to inform relay selections.

For the second case \cite{Sanaye-Pasand2009,Terzija2006,Rudez2016}, corresponding to adaptive UFLS techniques based on Wide-Area Measurement Systems (WAMS), a powerful monitoring and communications network is necessary. These WAMS-based adaptive UFLS methods use activation signals that are sent to relays in the instants following a contingency: real-time measurements from WAMS are typically needed for online estimations of the size of the contingency and level of inertia, in order to reduce the total amount of load to be shed.

Note that not all adaptive UFLS methods require WAMS, since adaptive UFLS refers to any UFLS method able to modulate the volume of load shed depending on the underlying system operating conditions at the time of activation, in order to achieve improved performance (e.g.~shedding less overall load). References \cite{localUFLS1,localUFLS2,localUFLS3,localUFLS4,localUFLS5,localUFLS6} proposed local yet adaptive UFLS methods, where WAMS is not required.

Since adaptive UFLS allocation is updated in real-time, it is unaffected by the uncertainty introduced by DG, as power flows through feeders could be available from the live measurements in WAMS. However, WAMS-based adaptive UFLS schemes rely completely on a communication network, which might fail or be subject to attacks \cite{Tondel2018}. Given that UFLS is a containment measure that avoids system collapse, robust methods to allocate feeders contributing to UFLS are needed, rather than purely relying on real-time updates. The method proposed in this paper can act as a back up to WAMS-based adaptive UFLS techniques, ensuring frequency security if the adaptive UFLS scheme cannot provide an acceptable performance due to unavailability of the communication network. 

Regarding the potential contribution of DG towards supporting frequency excursions, these devices have the capability to provide such support if operating with some headroom \cite{NREL_DER_FR}, and there is increasing interest in unlocking this functionality in real power grids. Other types of distributed resources such as electric vehicles also show promise in providing frequency containment services to the grid \cite{AliciaEV,MacEV}. Resorting to distributed resources for procuring grid services is however far from being an extended practice as of today, since virtually all system operators rely on transmission-level assets to provide most grid services, including frequency response. Nevertheless, if a certain feeder contained DG which provide frequency support, disconnecting that feeder in the UFLS scheme could not only imply a disconnection of net-generation from the grid, but also lose that valuable frequency support. Therefore, the system operator could simply remove these feeders from the initial list of candidate feeders fed into the UFLS-allocation algorithm, while remaining feeders (including those with DG not contributing towards frequency support), would be candidate feeders in the UFLS-allocation method introduced in the present work.

Given this background, the contributions of this work are:
\begin{enumerate}
    \item A convex-optimisation framework for the UFLS-allocation problem under uncertainty is proposed, enabling to meet a pre-defined level of risk in obtaining the minimum load to be shed, while accounting for correlation in net-load across the network.
    \item Two formulations for the chance constraint limiting risk are introduced, enabling to achieve the exact level of desired risk under Gaussian uncertainty, as well as a distributionally-robust approach, applicable to cases where uncertainty cannot be fully characterised due to lack of information, or uncertainty in different feeders is heterogeneous across the network.
    \item The advantages of the proposed method include providing the optimal selection of feeders for systems with traditional measurement-based UFLS schemes, as well as serving as a robust backup for WAMS-based adaptive UFLS schemes.
    \item Case studies on a test system are carried out to demonstrate the applicability of this feeder-selection procedure, as well as discussing the issue of social unfairness in the allocation process and proposing measures to mitigate it.
\end{enumerate}

The remainder of this paper is organised as follows: Section \ref{sec:Methodology} presents a description of the proposed methodologies for load shedding allocation under uncertainty. Section \ref{sec:CaseStudies} presents several case studies that demonstrate the applicability and advantages of said methodologies for UFLS allocation. Finally, Section~\ref{sec:Conclusions} gives the conclusion and suggests future lines of work.

\section{Methodology for risk-constrained UFLS allocation under uncertainty} \label{sec:Methodology}
This section formulates the task of selecting candidate feeders for the UFLS service, from a diverse set of feeders under net-power flow uncertainty, defined as a chance constrained mixed-integer optimisation problem. Two formulations to choose the optimal combination of feeders to meet the desired UFLS amount in a risk-secured manner are derived. Our presented method leverages distributional information on the uncertain feeders' net-power flow to do this in an efficient manner, 
whilst providing guarantees on a minimum volume of UFLS that would be delivered in the event of this service being triggered. 

Three formulations for the optimisation problem are discussed: 
\begin{enumerate}
    \item Deterministic allocation: represents the current status quo for selecting feeders that contribute toward the UFLS scheme, with the critical assumption that demand across the network follows fairly regular patterns. This assumption is becoming increasingly inaccurate as DG penetration increases. 
    \item Chance-constrained under Gaussian uncertainty: the first of our two proposed formulations, assumes that the net-power flows are Gaussian variables with known means and standard deviations. This is a restricting assumption for modelling uncertainty, though potentially applicable to some real systems, that would provide the widest possible feasible set for the chance-constrained optimisation.
    \item Distributionally robust chance-constrained: the second of our proposed formulations, assumes only knowledge of the net-power flows means and standard deviations, with no assumption on their specific probability distributions. Applicable to systems with diverse constituent feeders, where uncertainty arises from different sources such as DG and a variety of loads (e.g.~electric vehicles and heat pumps). The robustness of this method for accommodating any type of underlying uncertainty comes at the expense of some conservativeness in the optimisation.
\end{enumerate}

Each of these formulations is described in detail in following sections. Regarding the mathematical notation, decision variables in the optimisation problem are emboldened, while standard font is used for parameters. All vectors are defined as column vectors.

\subsection{Deterministic allocation approach} \label{sec:deterministic_approach}

Depending on the UFLS scheme deployed, the desired amount of load to be shed can be expressed as a set percentage of national demand or a set amount of demand. Here we formulate both problems mathematically and deterministically, before demonstrating their equivalence. This is relevant because it means the learnings from case studies performed for either scheme are transferable.

First, we consider the problem of selecting the optimal feeder combination to deliver a set percentage of national demand ($L^\%$). This is a simplistic deterministic method where the net-load through each feeder is assumed to always follow a certain percentage of national demand. This approach represents the current practice used extensively by system operators such as National Grid ESO in Great Britain \cite{WesternPower}. This formulation is defined by the following mixed-integer linear program:
\begin{subequations}
    \begin{alignat}{2}
    &\!\min_\mathbf{x}             &\quad&  (P^\%)^{T} \cdot \mathbf{x}  \\
    &\text{subject to } &     &   L^\% \leq (P^\%)^T \cdot \mathbf{x} \; \label{eq:deterministic_optimisation_problem_constraint_percentage} \\ 
    &                    &     & \mathbf{x}_i\; \in\; \{0,1\} \quad \forall \; i=1,\dots,m  
    \end{alignat}
    \label{eq:deterministic_optimisation_problem_percentage}
\end{subequations}
The objective is to minimise the percentage of load disconnected during the UFLS scheme, whilst ensuring that the minimum threshold for effective load shedding ($L^\%$), is met. Column vector `\mbox{$P^\%\; =\; [P^\%_1,P^\%_2,...,P^\%_m]^{T}\; \in \; \mathbb{R}^m$}' corresponds to the \% of national demand assumed for each of the $m$ participant feeders in the UFLS scheme. The penetration of DG within distribution networks make the forecast of this value uncertain. However, in this simple method that represents the current practice, the stochasticity is ignored and a fixed value is used (e.g.~the mean or a fixed percentile). Column vector `\mbox{$\mathbf{x}\; =\; [\mathbf{x}_1,\mathbf{x}_2,...,\mathbf{x}_m]^{T}\; \in \; \{0,1\}^m$}' corresponds to the binary decision variables in the optimisation (i.e.,~selection of a feeder contributing towards the UFLS scheme or not).

Second, the formulation for a set amount UFLS scheme is given by:
\begin{subequations}
    \begin{alignat}{2}
    &\!\min_\mathbf{x}             &\quad&  P^{T} \cdot \mathbf{x}  \label{eq:deterministic_Objective}\\
    &\text{subject to } &     &   L \leq P^T \cdot \mathbf{x} \; \label{eq:deterministic_optimisation_problem_constraint} \\ 
    &                    &     & \mathbf{x}_i\; \in\; \{0,1\} \quad \forall \; i=1,\dots,m  
    \end{alignat}
    \label{eq:deterministic_optimisation_problem}
\end{subequations}

This formulation is the same as in (\ref{eq:deterministic_optimisation_problem_percentage}) except that the column vector `$P\, =\, [P_1,P_2, ...,P_m]^{T}\, \\ \in \, \mathbb{R}^m$' now represents a set amount of UFLS expressed in MW, in the same way as the minimum UFLS requirement `$L$' is defined in MW. The reason for this equivalence is because translating a set amount of load to a set percentage of national demand is trivially achieved by dividing by the total national demand. Due to this equivalence, from here on only the set amount scheme will be considered in mathematical formulations and case studies. The importance and implications of moving to a MW-requirement for UFLS are explained in detail in Section~\ref{sec:Applicability}.

Before the rise of distributed RES, this deterministic method was appropriate for scheduling the UFLS services, as the net load through each feeder did follow fairly regular patterns. However, given the stochastic nature of renewable generation, this is no longer the case, and the likelihood of violating requirement (\ref{eq:deterministic_optimisation_problem_constraint}) could increase significantly. The key parameter in the deterministic formulation presented in this section is the forecast for net-load in each feeder, vector $P$. As will be shown in Section~\ref{sec:casestudy_deterministicapproach}, this formulation is not suitable for choosing a pre-defined risk in violating the minimum load-shedding threshold $L$: either an overly optimistic result (i.e.~lower total UFLS than $L$ would typically be available when this service is triggered) or overly conservative (i.e.~potentially much higher load being disconnected than the minimum amount required) would be achieved with this method. The UFLS-allocation formulations presented in following sections overcome this problem.

\subsection{Chance-constrained approach \label{sec:stochastic_approach}}

A chance-constrained optimisation approach is proposed here to account for the stochastic nature of feeders' net-load. The formulation in (\ref{eq:deterministic_optimisation_problem}) is modified by including a chance constraint which allows the system operator to specify an acceptable risk of UFLS under-delivery:
\begin{subequations}
    \begin{alignat}{2}
    &\!\min_\mathbf{x}            &\quad&  \mu_{\tilde{P}}^T \cdot \mathbf{x}  \\
    &\text{subject to } &     &  \mathbb{P}\left (   L \leq \tilde{P}^T \cdot \mathbf{x} \right ) \; \geq \; 1 - \epsilon \label{eq:cc} \\
    &                    &     & \mathbf{x}_i\; \in\; \{0,1\} \quad \forall \; i=1,\dots,m 
    \end{alignat}
    \label{eq:stochastic_optimisation}
\end{subequations}

The objective is to minimise the sum of expected load disconnection from all feeders, where column vector `\mbox{$\mathbf{\mu}_{\tilde{P}}\; =\; [\mu_1, \mu_2,...,\mu_m]^{T}\; \in \; \mathbb{R}^m$}' contains the mean values of net-load of each feeder participating in the UFLS scheme (i.e.~$\mathbf{\mu}_{\tilde{P}} = \textrm E[\tilde{P}]$). Constraint~(\ref{eq:cc}) ensures that the volume of UFLS disconnected is higher or equal to the minimum threshold `$L$', with a probability equal or greater than `\mbox{$ 1 - \epsilon $}'. In other words, `$\epsilon$' represents the risk of violating the minimum UFLS threshold. Column vector `\mbox{$\tilde{P} \; =\; [\tilde{P}_1,\tilde{P}_2,...,\tilde{P}_m]^{T}\; \in \; \mathbb{R}^m$}' contains the random variables for net-load in each feeder. 

The above formulation is non-deterministic and therefore impossible to solve.
However, if the means and covariances of $\tilde{P}$ are known (i.e.,~the first two moments of the random variables), moment-based distributionally robust chance constraint methods allow its reformulation into a tractable convex second-order cone problem \cite{Roald2015}. Here we present two methods depending on the distributional assumptions made on $\tilde{P}$.

We start by defining a new auxiliary variable:
\begin{equation} \label{delta def}
    \bm{\delta} = L - \tilde{P}^T \cdot \mathbf{x}
\end{equation}
This variable represents the cumulative under-delivery of UFLS compared to what is scheduled. It is desirable that this is negative with low probability, ensured when it is substituted into (\ref{eq:cc}).
\begin{equation} \label{delta cc}
    \mathbb{P}\left( \bm{\delta} \leq 0 \right) \; \geq \; 1 - \epsilon
\end{equation}
The mean and standard deviation of $\bm{\delta}$ are functions of the mean vector ($\mu_{\tilde{P}}$) and covariance matrix ($\Sigma$) of $\tilde{P}$:
\begin{equation} \label{delta mean}
    \mu_{\bm{\delta}} = L - \mu_{\tilde{P}}^T \cdot \mathbf{x}
\end{equation}
\begin{equation} \label{delta std}
    \sigma_{\bm{\delta}} = \sqrt{\mathbf{x}^T \Sigma \mathbf{x}} = \norm{ \Sigma^{\frac{1}{2}} \mathbf{x} }_2
\end{equation}

A covariance matrix is always symmetric and positive semi-definite. Accordingly, $\Sigma^{\frac{1}{2}}$ refers to the unique matrix that is positive semidefinite and such that $\Sigma^{\frac{1}{2}} \Sigma^{\frac{1}{2}} = (\Sigma^{\frac{1}{2}})^T \Sigma^{\frac{1}{2}} = \Sigma$.

We can normalise variable $\bm{\delta}$ to have zero mean and unity variance (leading to the new variable $\bm{\delta_n}$), and then substitute it into (\ref{delta cc}) to give:
\begin{equation} 
    \mathbb{P}\left( \bm{\delta_n} \leq \frac{-\mu_{\bm{\delta}}}{\sigma_{\bm{\delta}}} \right) \; \geq \; 1 - \epsilon
\end{equation}
The probability that $\bm{\delta_n}$ is less than or equal to a constant is given by `$F_\mathcal{P}(\cdot)$', the cumulative distribution function (CDF) of $\bm{\delta_n}$. Thus the above is equivalent to:
\begin{equation} \label{delta cc cdf}
     F_\mathcal{P} \left( \frac{-\mu_{\bm{\delta}}}{\sigma_{\bm{\delta}}} \right) \; \geq \; 1 - \epsilon 
\end{equation}
The exact form of the CDF depends on the assumed probability distribution of $\bm{\delta}$, which we refer to as `$\mathcal{P}$'. In this paper, we present the convex reformulation of (\ref{delta cc cdf}) under the assumption that the underlying uncertainty in net-power flow through feeders follows Gaussian distributions (i.e.,~$\mathcal{P}=\mathcal{N}(\mu_{\bm{\delta}},\sigma_{\bm{\delta}}^2))$, as well as a less specific distributionally robust assumption. Comparison of the efficiency and appropriateness of each of these assumptions is explored in the case studies.

\subsubsection{Convex reformulation of chance constraint under Gaussian uncertainty \label{sec:GA}}
Due to summative property of Gaussian distributions, if we assume that the net-power flow uncertainties ($\tilde{P}$) are Gaussian distributed, then so too will be $\bm{\delta_n}$. Thus, given the well defined CDF of Gaussian variables `$\Phi(\cdot)$', constraint  (\ref{delta cc cdf}) can be reformulated:
\begin{equation} \label{gauss 1}
     \frac{-\mu_{\bm{\delta}}}{\sigma_{\bm{\delta}}}  \; \geq \; \Phi^{-1}(1 - \epsilon) 
\end{equation}
Rearranging and substituting in the moment definitions:
\begin{equation} \label{gauss final form}
     \frac{\mu_{\tilde{P}}^T \cdot \mathbf{x} - L}{\Phi^{-1}(1 - \epsilon)}  \; \geq \; \sqrt{\mathbf{x}^T \Sigma \mathbf{x}} 
\end{equation}
Given that the under-delivery of UFLS is desired with low probability (certainly lower than 50\%), this equates to assuming `$\epsilon < 0.5$' and thus `$\Phi^{-1}(1 - \epsilon) > 0$'. Therefore, eq.~(\ref{gauss final form}) is a convex second-order cone (SOC) constraint of the standard form: 
\begin{equation}
    \norm{ A\mathbf{x} + b }_2 \leq c^T \mathbf{x} + d
\end{equation}
where: 
\begin{equation}
    A = \Sigma^{\frac{1}{2}} \quad 
    ; \quad 
    b = 0   \quad 
    ; \quad c = \frac{\mu_{\tilde{P}}}{\Phi^{-1}(1 - \epsilon)}  \quad 
    ; \quad 
    d = \frac{-L}{\Phi^{-1}(1 - \epsilon)}
\end{equation}

The derived constraint compels the optimiser to find a combination of feeders that will result in at least `$L$' MW of load shed with `$(1-\epsilon)$\%' probability. If a given feeder `$i$' is chosen by the optimal allocation, making its corresponding binary variable equal unity ($\mathbf{x}_i=1$), its mean value will increase the left-hand side of (\ref{gauss final form}), whilst its variance and covariance coefficients will increase the right-hand side. Thus, the mathematical structure of constraint (\ref{gauss final form}) is aligned with the intuitive result that feeders with a high mean and low uncertainty (i.e.~low DG penetration) are most desirable.

%%%%%%%%%%%%%%%%%%% CC DISTRIBUTIONALLY-ROBUST %%%%%%%%%%%%%%%%%%%

\subsubsection{Convex formulation of a distributionally-robust chance constraint \label{sec:DR}}

The expression in~(\ref{gauss final form}) describes the reformulation of the chance constraint (\ref{eq:cc}) when the probability distribution of $\tilde{P}$, and therefore $\bm{\delta}$, is Gaussian. However, this is a restrictive assumption that in some situations need not hold. In this section, a second reformulation of the chance constraint is derived. In this case, the objective is to have an expression that guarantees compliance with the constraint when limited information about the characterisation of uncertainty in net-load in each feeder ($\tilde{P}$) is available. We assume here that only the mean and standard deviation of each component in $\tilde{P}$ is known, but no information on the PDFs is available. The  formulation presented in this section therefore provides a convex reformulation for the chance constraint for cases when at least some of the feeders' net-load does not follow Gaussian distributions, or information is lacking on their actual distribution.

According to (\ref{delta mean}) and (\ref{delta std}), knowledge of the mean and standard deviation of the net-load in each feeder implies knowledge of the mean and standard deviation of $\bm{\delta}$. However, the lack of knowledge on the distributions of $\tilde{P}_i$ makes the true distribution of $\bm{\delta}$ ambiguous. In other words, the true distribution `$\mathcal{P}$' of $\bm{\delta}$ could be any within the set of distributions with $\mu = \mu_{\bm{\delta}}$ and $\sigma = \sigma_{\bm{\delta}}$, called the ambiguity set ($\mathcal{A}$). 

This ambiguity inhibits the use of an exact form of $\bm{\delta}$'s CDF being used to reformulate (\ref{delta cc cdf}), as was done for the Gaussian case in Section~\ref{sec:GA}. Instead, (\ref{delta cc}) should be reformulated into convex form such that it always holds, even for the worst case distribution within $\mathcal{A}$. In other words, a distributionally robust formulation is sought.

Defining a lower bound on `$F_\mathcal{P}(\cdot)$', the CDF of $\mathcal{P}$ within $\mathcal{A}$, allows the inequality within (\ref{delta cc}) to be maintained and progress to be made. For some positive constant $\lambda$:
\begin{equation} \label{cdf lower bound}
  f_\mathcal{P}(\lambda) = \textup{inf}_{\mathcal{P} \in \mathcal{A}} \ F_\mathcal{P}(\lambda)
\end{equation}
\begin{equation}  \label{cdf equivalence}
    F_\mathcal{P}(\lambda) \geq  f_\mathcal{P}(\lambda) \geq 1-\epsilon
\end{equation}
Following the method presented in \cite{Roald2015}, the classical Cantelli inequality can be used to find the appropriate form of $f_\mathcal{P}(\lambda)$:
\begin{equation}  \label{cdf cantelli}
    f_\mathcal{P}(\lambda) = \frac{\lambda^2}{1 + \lambda^2}
\end{equation}
With its corresponding well-defined inverse function, given that a CDF is monotonically increasing:
\begin{equation} \label{inverse cdf}
    f_\mathcal{P}^{-1}(\lambda) = \sqrt{\frac{\lambda}{1-\lambda}}
\end{equation}
Thus by using (\ref{inverse cdf}) in (\ref{cdf equivalence}) with `$\lambda = 1-\epsilon$' and rearranging, the distributionally robust form of the chance constraint (\ref{eq:cc}) is found:
\begin{equation} \label{DRO final form}
     \sqrt{\frac{\epsilon}{1-\epsilon}} \, \cdotp (\mu^T\mathbf{x} - L)  \; \geq \; \sqrt{\mathbf{x}^T \Sigma \mathbf{x}} 
\end{equation}

Comparing between equations~(\ref{gauss final form}) and~(\ref{DRO final form}) reveals that they are the same, except for the constant on the left-hand side, thus (\ref{DRO final form}) is also a convex SOC with wide applicability. Some expectations of the performance of the constraints can be deduced from their mathematical formulation. For any $\epsilon$:

\begin{equation}
    \sqrt{\frac{1-\epsilon}{\epsilon}} > \Phi^{-1}(1-\epsilon)
    \label{eq:dr_bigger}
\end{equation}

\begin{figure}[!t]
    \centering
    \includegraphics[width=0.85\textwidth]{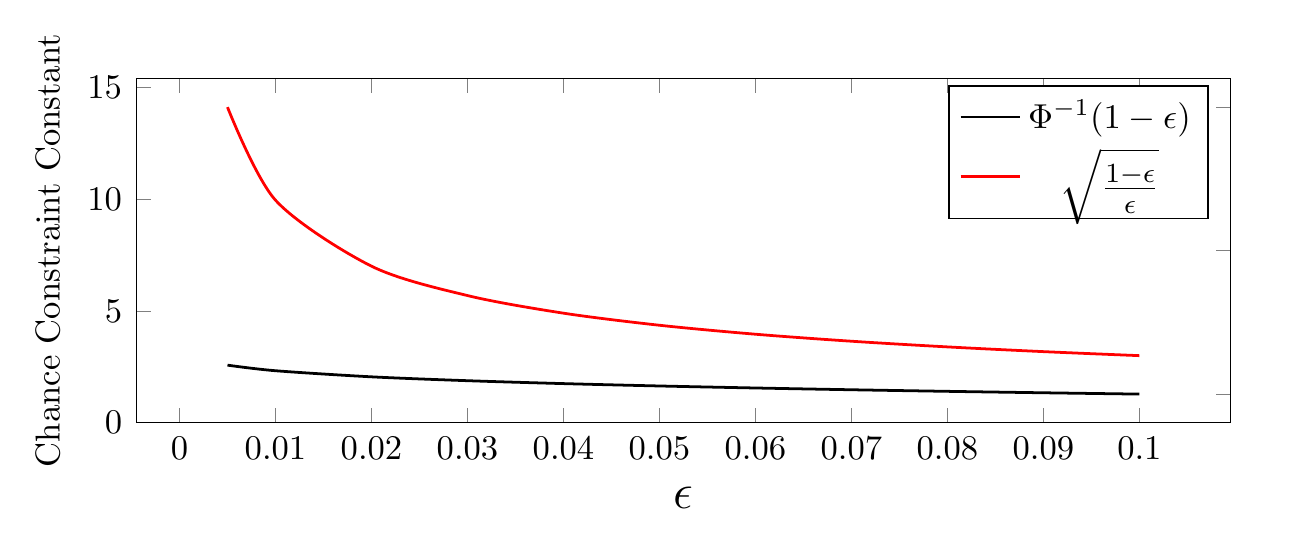}
    \caption{Value of the left-hand side constant for chance constraints (\ref{gauss final form}) and  (\ref{DRO final form}), over a range of risk preference ($\epsilon$). To accommodate any distribution type with a specified mean and standard deviation, the distributionally robust constraint's constant is larger over all $\epsilon$, resulting in a robust and therefore more conservative solution.}
    \label{fig:epsilon}
\end{figure}

This means that for the distributionally-robust constraint, `$\mathbf{\mu^T} \cdot \mathbf{x} - L$' must be a greater quantity than in the Gaussian assumption chance constraint. Thus, the distributionally-robust constraint yields results that are more conservative, i.e.~more load could potentially be disconnected in order to assure meeting the required reliability. Figure~\ref{fig:epsilon} plots the relationship between these two constants and shows that the relative conservativeness of  (\ref{DRO final form}) over (\ref{gauss final form}) increases as the operator seeks higher delivery confidence (i.e. smaller $\epsilon$).

\subsection{Applicability to current and future scenarios of power system operation} \label{sec:Applicability}

\begin{figure}[!t]
\vspace*{-1mm}
    \centering
    \includegraphics[width=0.88\textwidth]{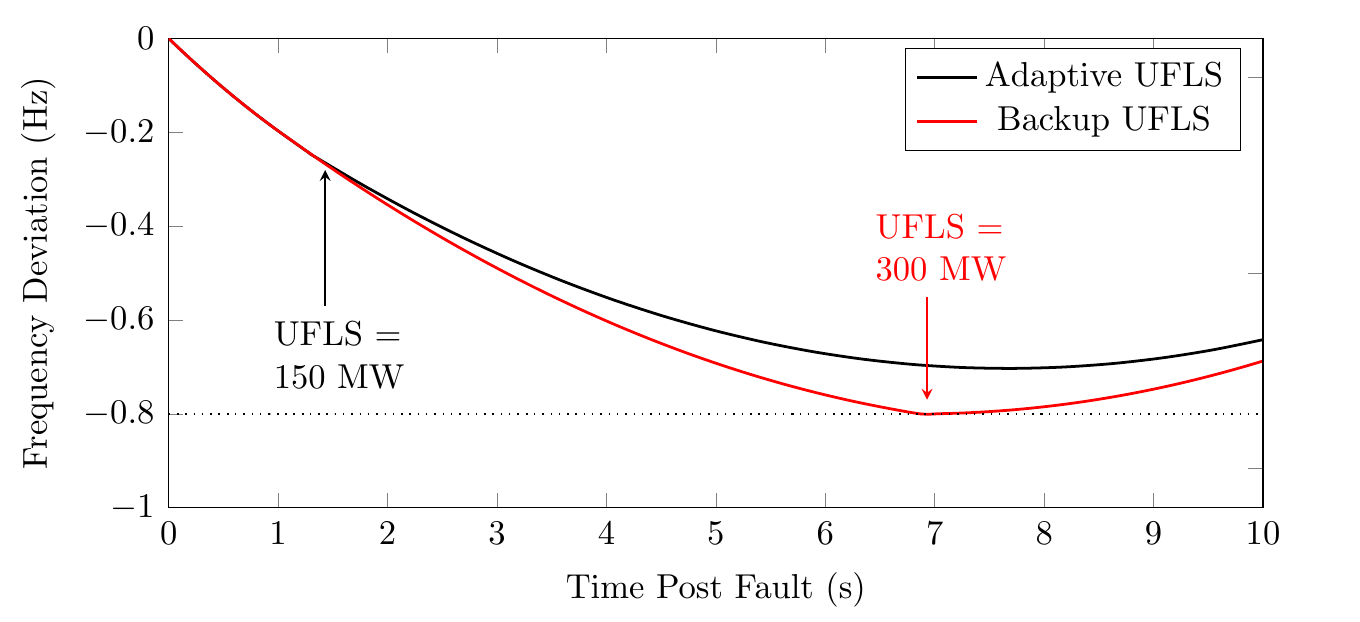}
    \caption{Frequency excursions following a generation loss, for a WAMS-based adaptive UFLS scheme with operational and unavailable communications. When communications fail, the relays selected for the backup scheme disconnect their load once their local frequency reaches the 0.8Hz trigger level. These relays must be selected ahead of time and thus are subject to net-load uncertainty caused by DG.}
    \label{fig:dynamics}
\end{figure}

The methods presented above facilitate the selection of the optimal feeder combination to provide the required UFLS amount, given uncertainty in net-power flow through feeders. We now discuss the applicability of these methods to real power system operation, focusing on the need for a communication network to update UFLS-relay settings with certain periodicity, as well as the role that the proposed methods would play, either as primary or secondary set-point for the UFLS relays. For this, it is relevant to consider different current and future scenarios for power system operation:

\begin{enumerate}
   \item Current practice for UFLS, where relays are tuned very infrequently and the goal is to shed a set percentage of national demand in the event of a large frequency excursion \cite{Rowland}: an uncertainty-aware method for selecting the feeders contributing towards the UFLS scheme, such as the one proposed in this paper, will be necessary as DG penetration increases. This operational practice has traditionally relied on the assumption that electricity consumption is roughly uniform across the network, which will become increasingly inaccurate with the integration of DG. However, this operational practice could continue if combined with a risk-aware UFLS allocation such as the one discussed here, although it could potentially lead to very conservative outcomes: given that the uncertainty in net-load would have to be considered with a very long look-ahead (e.g.~one year, the periodicity for re-tuning UFLS relays), the variance could be very large, leading to very conservative UFLS-allocation results (i.e.~a large amount of load being disconnected if the UFLS service is triggered).
   
   Given that no communication network is necessary for this operational practice, relays could be updated manually and very infrequently (e.g.~once every year). Then, the method introduced in this paper would find the best combination of feeders to be chosen for the UFLS scheme, given the uncertainty in their instantaneous power consumption for the next year.

   \item `Fully controllable' future operation, where a powerful control and communication network has been deployed and therefore adaptive UFLS is widely available: the method proposed in this paper could be used as a backup set-point for UFLS relays. Given that WAMS-based adaptive UFLS methods rely on demanding communication requirements, since the set-point of UFLS relays must be updated in a timescale of 1 or 2 seconds following a contingency, it would be highly valuable to make these methods robust against failures in the communication network. The backup set-point determined by the uncertainty-aware UFLS-allocation method we propose here would be based on local frequency measurements (instead of control signals like WAMS-based adaptive UFLS), that are the last resort if the communication network needed for the adaptive UFLS scheme has failed. 
   
   The usefulness of our proposed method as a backup to a WAMS-based adaptive UFLS scheme is demonstrated in Figure \ref{fig:dynamics}, which plots an adaptive UFLS scheme with and without successful communication. When the communication network is fully operational, the controller detects the fault and trips 150~MW of load at 1.4~s. This shedding of load combines with system inertia and frequency response services to arrest frequency decline within the 0.8~Hz limit that would trigger conventional UFLS relays. However, if communication fails, no load is tripped at 1.4~s and thus frequency declines to 0.8~Hz. Relays detect this via local frequency measurements and trip 300~MW of load to restore demand-generation balance and arrest frequency decline. The constituent feeders for this `backup' UFLS must be chosen ahead of time and are subject to the same net-load uncertainty from DG as traditional UFLS schemes, therefore a risk-aware feeder selection method such as the one proposed in this paper would be needed.

   \item `Frequently tunable' future operation, where a communication network is in place, but the system operator does not rely on updating UFLS-relay settings in a sub-second timescale (i.e.~WAMS-based adaptive UFLS is not adopted): our feeder-allocation method would enable optimal operation of the UFLS scheme in this scenario. The communication network would be used to frequently update the settings of UFLS relays (e.g.~every hour), where the risk-aware UFLS-allocation method would be used to select the appropriate feeders, given the uncertainty in net-load through each feeder in the next hour. Based on this set-point updated hourly, local frequency measurements would make the UFLS relays trip if the frequency excursion exceeds the pre-defined threshold. As will be demonstrated through relevant case studies in Section~\ref{sec:CaseStudies}, this method provides risk guarantees by leveraging the available information on the uncertainty in net-load.
\end{enumerate}
Furthermore, the chance-constrained allocation method would provide optimal performance in the above scenarios 2 and 3,
when working in tandem with a UFLS-scheduling method such as the one proposed in \cite{Omalley}. First, the optimal amount of UFLS to be scheduled would be found following the procedure detailed in \cite{Omalley}, and this amount of UFLS would then be optimally allocated among feeders using the chance-constrained method proposed in this paper. Given that the triggering mechanism for UFLS relays would be local measurements (i.e.~the relays trip when the local frequency drops below the threshold), the communication network would only be used immediately following scheduling to select the feeders (e.g.~every hour or half-hour, each time the energy and ancillary services markets are cleared), and thus this procedure is robust to communication delays or failure during the safety-critical situations that require UFLS. Nevertheless, the updating of UFLS relays could also be done sub-hourly, if information on the net-load through feeders is received with sub-hourly resolution: the closer to real-time the selection of feeders is updated, the lower the forecast error for DG power output will be, therefore the lower the conservativeness of the chance-constrained UFLS allocation will be.

Finally, regarding the applicability of this UFLS-allocation method to a large power system, such as the Continental Europe synchronous area, it is worth considering the case where different system operators in the same synchronous areas apply slightly different criteria for their UFLS requirements. Furthermore, after a large disturbance, the synchronous area could split into different electric islands due to loss of synchronism, and the UFLS scheme would need to act as the last resort for containing the generation-demand imbalance within each island, even though the level of inertia in each island would not be known at the time of selecting UFLS feeders. In this context, the UFLS-allocation method proposed in this paper would still perform appropriately: the volume of load shed would likely be higher than optimal, given that the optimal amount can only be computed for a known value of inertia and power imbalance, as explained in \cite{Omalley}. The feeders selected by this algorithm applied by each system operator following their own criteria, such as specific trigger levels for UFLS, would be disconnected in a low-inertia island once frequency drops below each of the UFLS trigger levels. Therefore, it is guaranteed that sufficient load would be disconnected, even if more load than absolutely necessary would be shed, given that the inertia available and power imbalance in the island was not known a priori.

\section{Case studies} \label{sec:CaseStudies}

In this section, the performance of the methodologies for load shedding allocation introduced in Section~\ref{sec:Methodology} is analysed for a test system with 20 candidate feeders with varying degrees of uncertainty, described in Table \ref{tab:feedersdata}. 
The case studies were run with a minimum permissible load shed $L$ of $250\textrm{MW}$. It should be noted that because the mathematical formulation for the `set percentage' UFLS scheme is the same as the `set amount' scheme (as shown in Section~\ref{sec:deterministic_approach}), the learnings here are transferable.  All optimisation problems in this section were modelled using the YALMIP toolbox for MATLAB \cite{LOFBERG2004}, calling Gurobi as the numerical solver while setting the acceptable integer optimality gap to $0.01\%$.

\begin{figure}[!htb]
\minipage{0.5\textwidth}
  \includegraphics[width=1.13\linewidth]{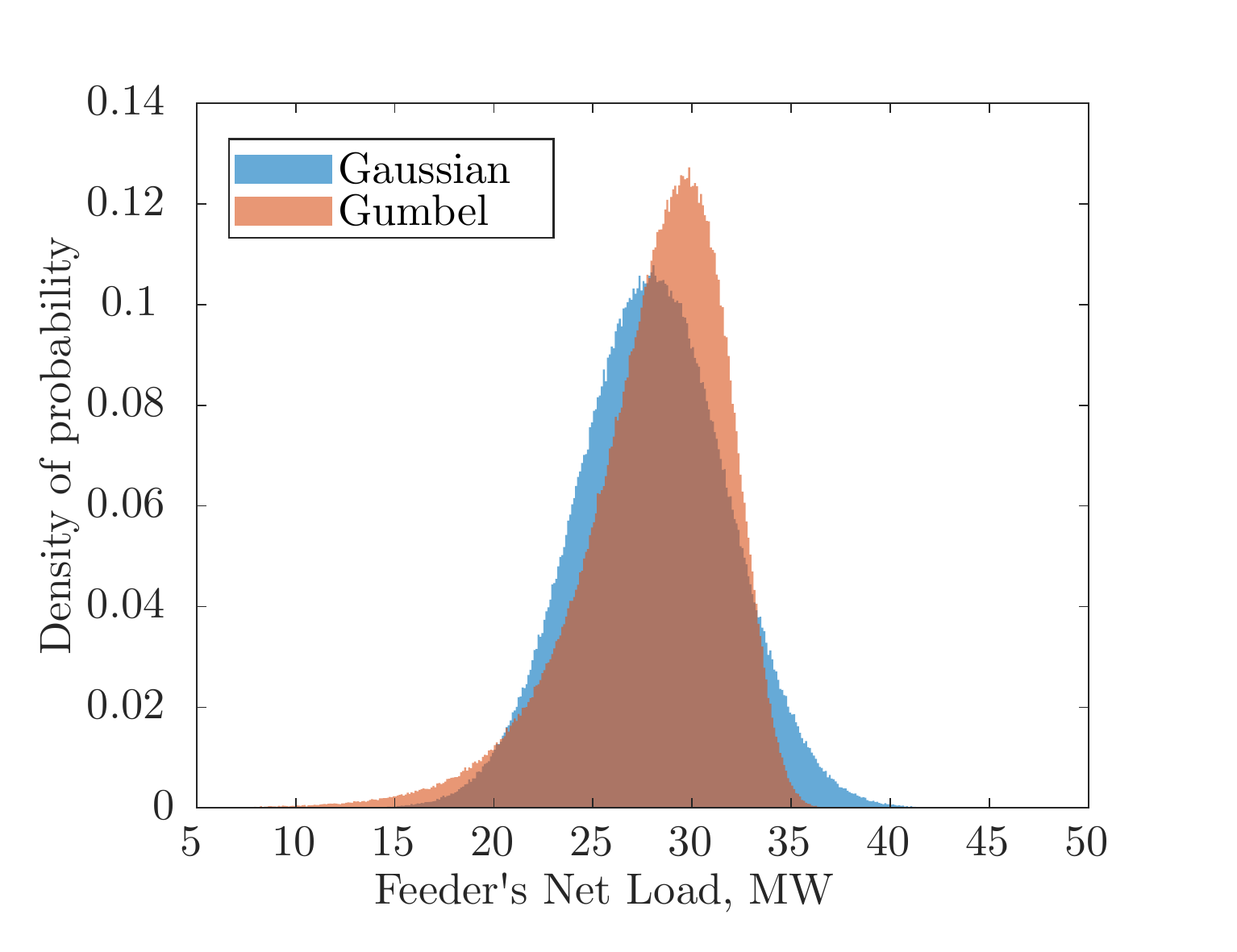}
\endminipage\hfill
\minipage{0.5\textwidth}
  \includegraphics[width=1.13\linewidth]{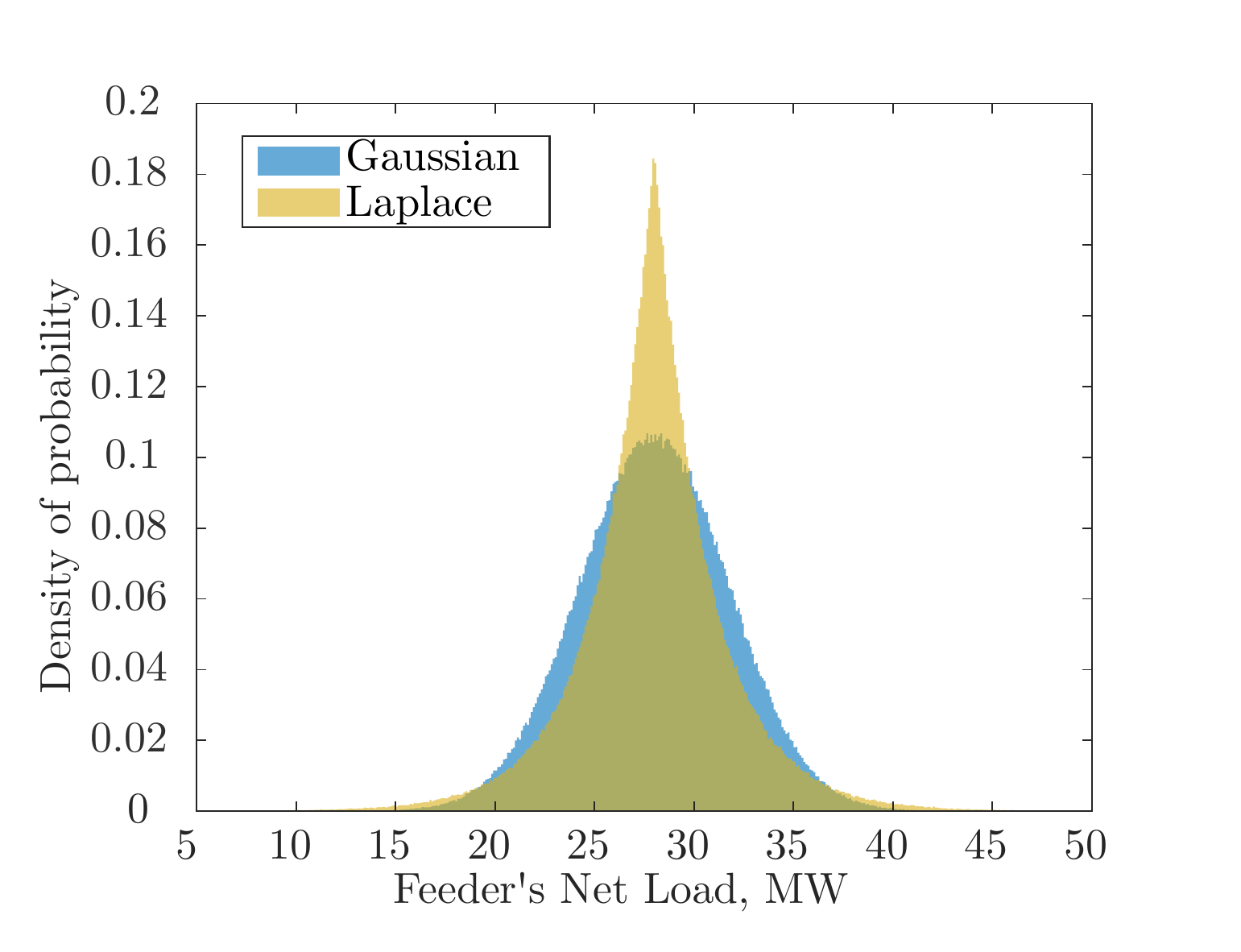}
\endminipage\hfill
\centering
\minipage{0.5\textwidth}
  \includegraphics[width=1.13\linewidth]{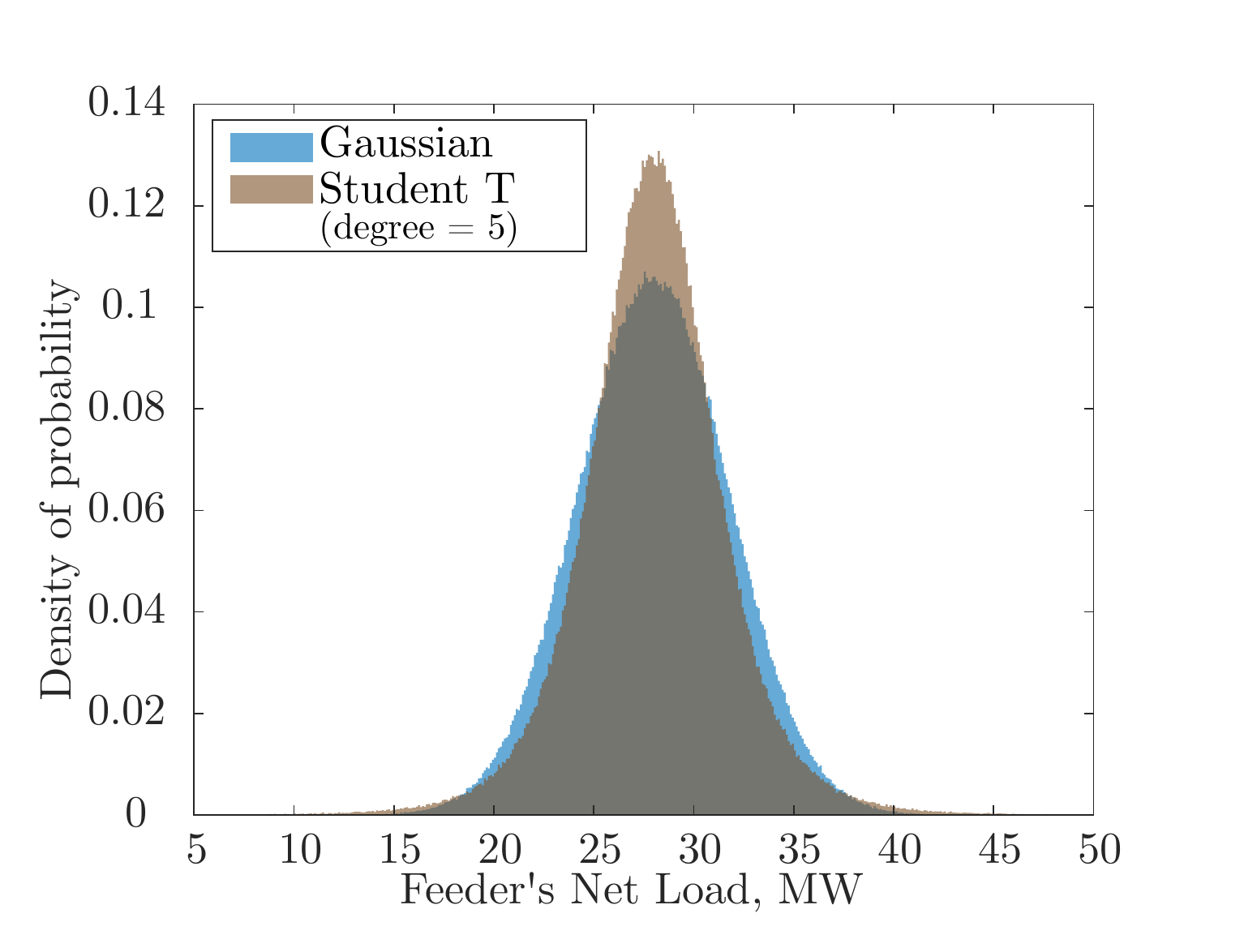}
\endminipage
  \caption{Comparison between the Gaussian and Gumbel, Laplace and Student's \textit{t} distributions, all with the same values for the first two moments.}
  \label{fig:comparison_distributions}
\end{figure}

We consider different underlying probability distributions for the uncertainty in feeder's net-load. First, the Gaussian PDF is taken as a reference due to its common use in probabilistic forecasting of both load and DG in power systems \cite{Hasan2019,arilezasoroudi2014,ZAKARIA20201543}. Second, the Gumbel distribution
is used to consider higher probability of low net-load values (i.e.,~greater DG output), which has been applied to modelling  uncertainty in electric load and renewable generation \cite{Hasan2019,LI2020116699}. Finally, the Laplace and Student's \textit{t} distributions were chosen to model heavy-tailed probabilistic forecasts, as these distributions are frequently used in the literature to model uncertainty in solar and wind generation \cite{DONG2020333,WEIJUN2018}. Figure~\ref{fig:comparison_distributions} shows a comparison of these probability distribution functions with the same value for the first two moments (i.e.~mean and standard deviation), of $\mu = 27.5 \textrm{MW}$ and $\sigma = 3.75 \textrm{MW}$. The mathematical characterisation of the PDFs is given in Table~\ref{tab:distirbution_parameters}, where the location and scale parameters were calculated as a function of the mean and standard deviation.  

\begin{table}[H]
\centering
\small
    \caption{Test system with 20 candidate feeders for the UFLS scheme, with varying levels of uncertainty in net-load of each feeder.}
    \begin{tabular}{c|c|c}
    \multicolumn{1}{l|}{\textbf{Feeder \#}} & \multicolumn{1}{l|}{\textbf{$\mu_{\tilde{P}}$ (MW)}} & \multicolumn{1}{l}{\textbf{$\sigma_{\tilde{P}}$ (MW)}} \\ \hline
    1 & 10.00 & 4.34 \\
    2 & 21.00 & 3.06 \\ 
    3 & 34.00 & 4.88 \\ 
    4 & 30.00 & 3.21 \\ 
    5 & 16.00 & 2.62 \\ 
    6 & 35.00 & 3.44 \\ 
    7 & 33.00 & 3.41 \\ 
    8 & 20.00 & 3.92 \\ 
    9 & 29.00 & 1.04 \\ 
    10 & 13.00 & 4.42 \\ 
    11 & 22.00 & 2.10 \\ 
    12 & 38.00 & 2.28 \\ 
    13 & 28.00 & 3.48 \\ 
    14 & 22.00 & 4.95 \\ 
    15 & 27.00 & 4.00 \\ 
    16 & 40.00 & 4.51 \\ 
    17 & 34.00 & 3.72 \\ 
    18 & 19.00 & 3.41 \\ 
    19 & 18.00 & 1.51 \\ 
    20 & 16.00 & 1.77 \\ \hline
    \end{tabular}
\label{tab:feedersdata}
\end{table}

\begin{table}[H]
    \caption{Mathematical description of the different PDFs considered.}
    \begin{tabular}{ccll}
    \multicolumn{1}{c}{} & \multicolumn{1}{c}{\multirow{2}{*}{\textbf{\begin{tabular}[c]{@{}c@{}}Probability \\ Distribution Function\end{tabular}}}} & \multicolumn{2}{c}{\textbf{Parameters}} \\ \cline{3-4} 
    \multicolumn{1}{c}{} & \multicolumn{1}{c}{} & \multicolumn{1}{c}{\textbf{Location}} & \multicolumn{1}{c}{\textbf{Scale}} \\ \hline
    \multicolumn{1}{c}{\textbf{\begin{tabular}[c]{@{}c@{}}Gaussian\end{tabular}}} & \multicolumn{1}{l}{$f(x | \alpha , \beta) = \frac{1}{\sigma \cdot \sqrt{2\pi} } \cdot e^{-\left ( \frac{x- \mu}{\sigma}  \right )^2 }$} & \multicolumn{1}{l}{$\alpha = \mu$} & \multicolumn{1}{l}{$\beta = \sigma$} \\ \hline 
    \multicolumn{1}{c}{\textbf{\begin{tabular}[c]{@{}c@{}}Gumbel\end{tabular}}} & \multicolumn{1}{l}{$f(x | \alpha , \beta) = \frac{1}{\beta } \cdot e^{\frac{x- \mu}{\beta} - e^{\frac{x-\mu}{\beta} } }$} & \multicolumn{1}{l}{$\alpha = \mu + \beta \cdot \gamma^{*} $} & \multicolumn{1}{l}{$\beta = \frac{\sqrt{6}}{\pi} \cdot \sigma$} \\ \hline
    \multicolumn{1}{c}{\textbf{\begin{tabular}[c]{@{}c@{}}Laplace\end{tabular}}} & \multicolumn{1}{l}{$f(x | \alpha, \beta) = \frac{1}{2 \cdot \beta } \cdot e^{\frac{|x- \mu|}{\beta}}$} & \multicolumn{1}{l}{$\alpha = \mu$} & \multicolumn{1}{l}{$\beta = \frac{\sigma}{\sqrt{2}}$} \\ \hline
    \multicolumn{1}{c}{\textbf{\begin{tabular}[c]{@{}c@{}}Student's \textit{t}\end{tabular}}} & \multicolumn{1}{l}{$f(x | \nu ^{**} ) = \frac{\Gamma \left ( \frac{\nu+1}{2} \right )}{\sqrt{\nu \cdot \pi} \cdot \Gamma \left ( \frac{\nu}{2} \right )} \cdot \left ( 1 + \frac{x^2}{\nu} \right )^{-\frac{\nu+1}{2}}$} & \multicolumn{1}{l}{\begin{tabular}[c]{@{}l@{}}$\alpha = 0,$\\ $\textrm{for} \; \nu>0$\end{tabular}} & \multicolumn{1}{l}{\begin{tabular}[c]{@{}l@{}}$\beta = \frac{\nu}{\nu-2},$\\ $ \textrm{for} \; \nu > 2 $\end{tabular}} \\ \hline
    \multicolumn{4}{l}{\small\begin{tabular}[c]{@{}l@{}}$^{*}$ Euler-Mascheroni constant\\ $^{**}$ Degree of freedom\end{tabular}}
    \end{tabular}
    \label{tab:distirbution_parameters}
\end{table}
\vspace*{-1mm}

\subsection{Limitations of a deterministic optimisation approach} \label{sec:casestudy_deterministicapproach}

First we consider the deterministic optimisation approach defined in (\ref{eq:deterministic_optimisation_problem}) for allocating UFLS under uncertainty, as described in Section~\ref{sec:deterministic_approach}. For this approach, a fixed percentile of the net-load forecast in each feeder must be chosen before solving the optimisation, as this formulation does not enable the operator to pre-specify the desired level of risk for the overall UFLS-allocation result. That is, if the system operator is willing to accept a 1\% chance of violating the minimum UFLS threshold `$L$', using this deterministic formulation implies that a guess must be made on the percentile of net-load in each feeder that would lead to the desired level of overall risk.

To illustrate this point, six different percentiles of net-load for each feeder described in Table~\ref{tab:feedersdata} are found and input into (\ref{eq:deterministic_optimisation_problem}). Without loss of generality and for the sake of simplicity, Gaussian uncertainty for the forecast of net-load is assumed for now. The $1\%$, $10\%$, $20\%$, $30\%$, $40\%$ and $50\%$ percentiles of this PDF are found for each of the 20 feeders, and used as the values within vector $P$ in eqs.~(\ref{eq:deterministic_Objective}) and (\ref{eq:deterministic_optimisation_problem_constraint}). The results of solving the optimisation problem in (\ref{eq:deterministic_optimisation_problem}) with each of these percentiles are included in Table~\ref{tab:results_deterministic}, where two metrics are shown: 1) the percentage of occurrences that would violate the load disconnection requirement `$L$', computed by sampling 100,000 times from the Gaussian PDF for each feeder chosen in the optimisation; and 2) the expected load disconnection, which corresponds to the sample mean of the sampling process carried out to obtain the previous metric.

\begin{table}[h!]
\centering
\small
\caption{Results of allocation of load shedding using a deterministic optimisation approach.}
    \begin{tabular}{c|c|c}
    \begin{tabular}[c]{@{}c@{}}Optimisation \\ percentile (\%)\end{tabular} & \begin{tabular}[c]{@{}c@{}}Samples below \\ disconnection \\ requirement (\%)\end{tabular} & \begin{tabular}[c]{@{}c@{}}Expected  load \\ disconnection (MW)\end{tabular} \\ \hline
    1 & 0.00 &  374.02 \\ 
    10 & 0.01 &  289.01 \\ 
    20 & 0.30 &  282.97 \\ 
    30 & 5.89 &  265.96 \\ 
    40 & 23.71 &  256.96 \\ 
    50 & 50.08 &  249.98 \\ \hline
    \end{tabular}
    \label{tab:results_deterministic}
\end{table}

The results demonstrate that using a percentile of 1\% for each feeder leads to an extremely conservative allocation of load shedding, since the combined probability of all feeders falling below the $1^\textrm{th}$ percentile forecast error for net-load is almost negligible. To achieve an overall risk of 1\%, a percentile between 20\% and 30\% for each feeder must be chosen. This highlights the unsuitability of the deterministic scheme when dealing with stochastic net-load through feeders in future systems with high DG penetration: a system operator is only concerned with the aggregate load shed being guaranteed above a certain level, but there is no fixed relationship between the individual percentile chosen and this aggregate security level. This means the operator can only guarantee security by being overly conservative and incurring very high levels of potential load shed.

Therefore, a deterministic optimisation approach, which is fairly representative of current practice in most systems, is inadequate for the UFLS-allocation problem under significant DG penetration. Future power systems will require to explicitly account for uncertainty within the UFLS-allocation problem.

\subsection{Chance-constrained approach}

\subsubsection{Impact of correlation in net-load probabilistic forecasts across the network}\label{sec:results_correlation}

This section demonstrates the performance of the chance-constrained approach in the allocation of load shedding, which not only allows to achieve a pre-defined level of risk, but notably enables consideration of correlation in net-load of the different feeders. Since net-load of active feeders is largely driven by DG output, which in turn is driven by meteorological conditions, correlation in feeders that are geographically close is an important characteristic that must be accounted for \cite{ZHANG2018743}. Again, Gaussian uncertainty is used here for simplicity, while the learnings on the importance of considering correlations are applicable to any other probability distributions.

\begin{figure}[H]
\minipage{0.5\textwidth}
  \includegraphics[width=1\linewidth]{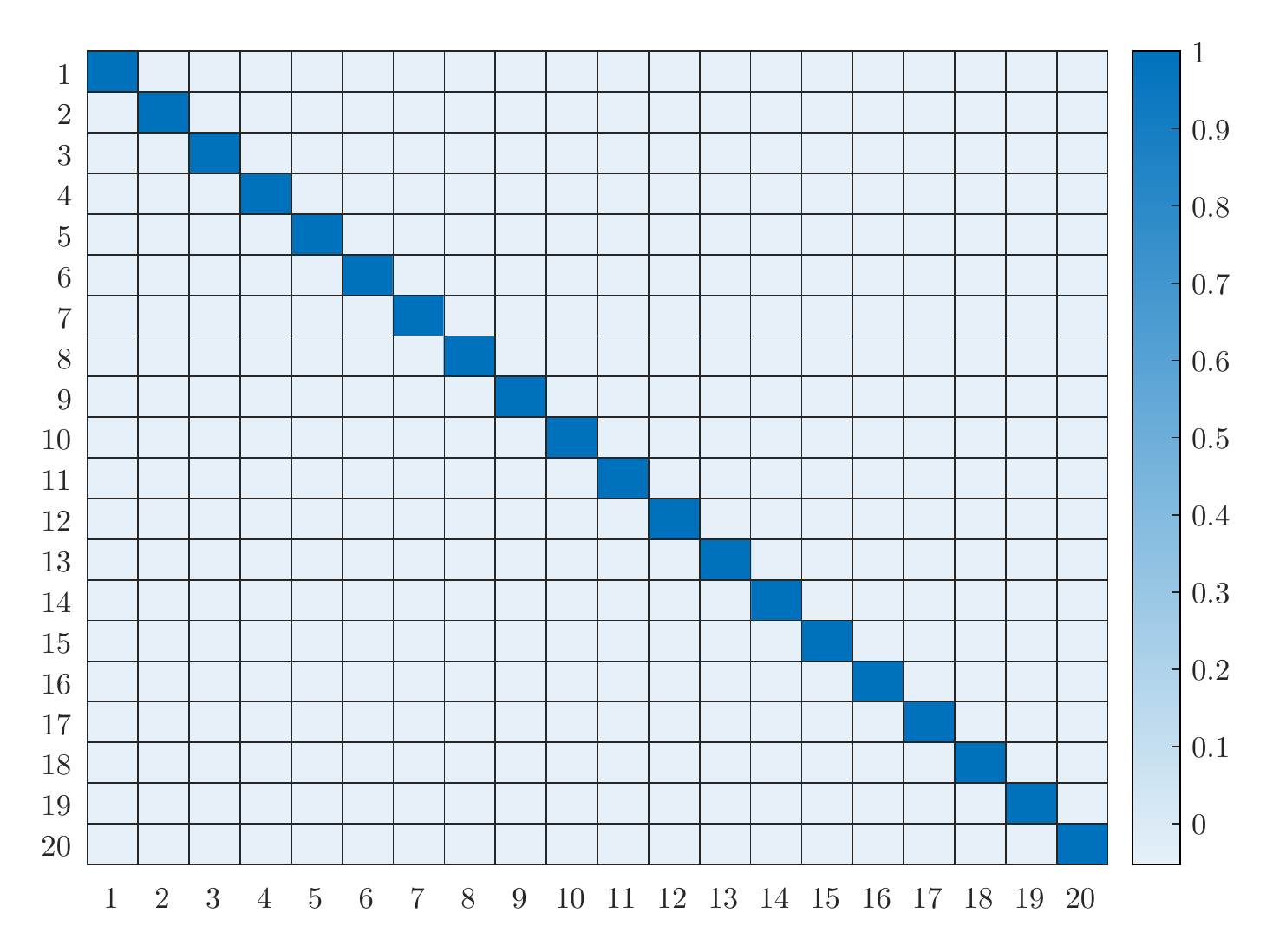}
\endminipage\hfill
\minipage{0.5\textwidth}
  \includegraphics[width=1\linewidth]{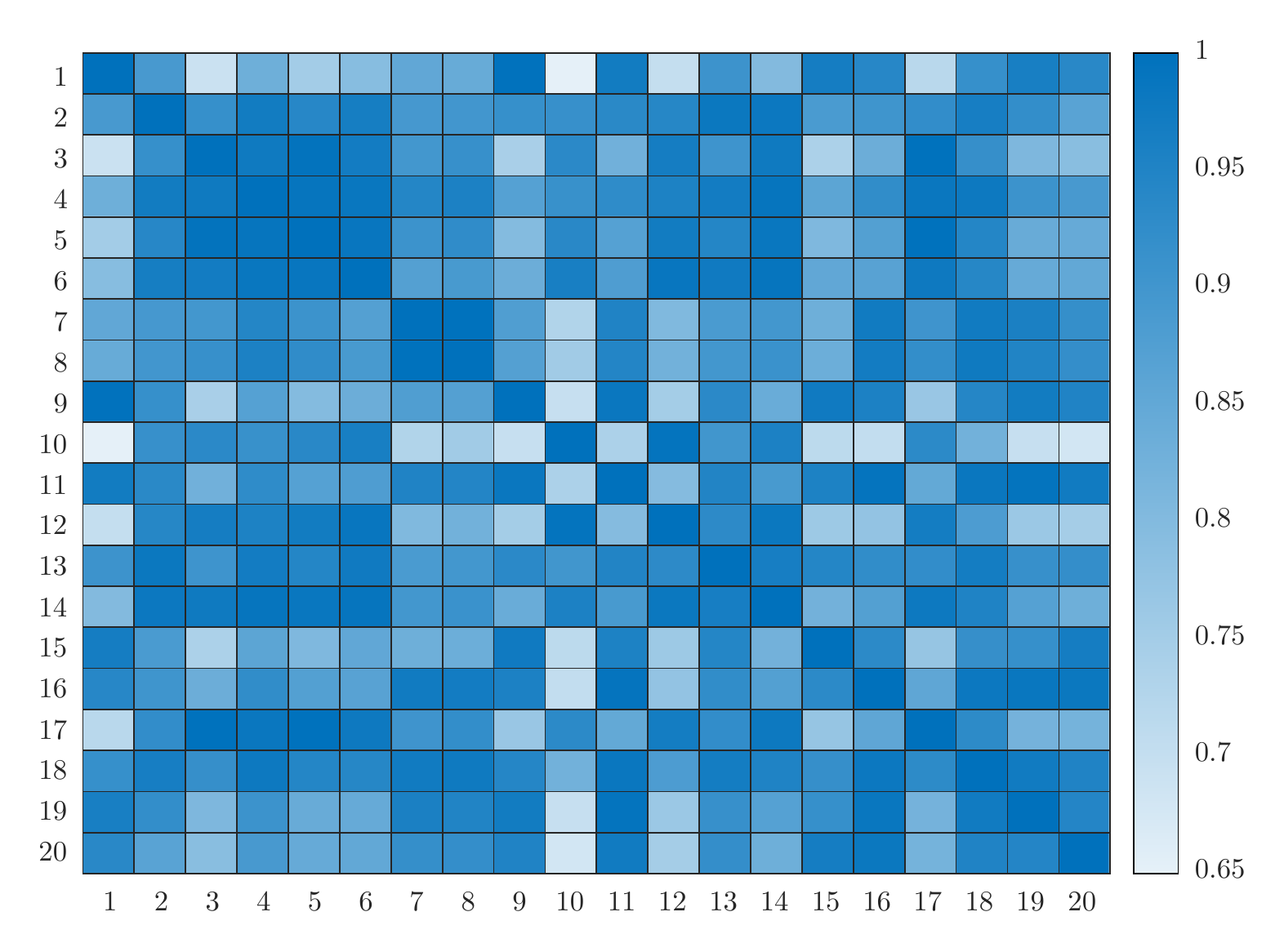}
\endminipage\hfill
  \caption{Heatmap for the covariance matrices of the uncorrelated and correlated net-load probabilistic forecast.}
  \label{fig:correlation_comparison}
\end{figure}

To demonstrate the importance of using appropriate probabilistic forecasts in the UFLS-allocation method, that account for correlation, here we solve the chance-constrained optimisation defined in~(\ref{eq:stochastic_optimisation}), while considering two different covariance matrices `$\Sigma$', in eq.~(\ref{gauss final form}). Heatmaps for these covariance matrices are depicted in Figure \ref{fig:correlation_comparison}. 
Table~\ref{tab:results_correlation} shows the results for two cases: Case 1 assumes uncorrelated random variables for the net-load of feeders (i.e.,~the covariance matrix is diagonal, containing only the variance of each feeder's net-load), while Case 2 does consider correlation in the optimisation problem. The sampling for feeders' net-load performed after solving the optimisation does account for the correlation in both cases, since DG outputs are in fact correlated. The results show that when correlation is neglected as in Case 1, the optimisation is overly optimistic: the UFLS-threshold defined in constraint~(\ref{gauss final form}) is met for only $19\%$ of the samples, while the accepted risk was set to 1\%.

\begin{table}[h!]
\small
\centering
    \caption{Results of neglecting (Case 1) and considering (Case 2) correlation among feeders' net-load forecasts.}
    \begin{tabular}{lccc}
    \multicolumn{1}{c|}{} & \multicolumn{1}{c|}{$\epsilon$   (\%)} & \multicolumn{1}{c|}{\begin{tabular}[c]{@{}c@{}}Samples below\\ disconnection\\ requirement (\%)\end{tabular}} &  \multicolumn{1}{c}{\begin{tabular}[c]{@{}c@{}}Expected load \\  disconnection (MW)\end{tabular}} \\ \hline
    \multicolumn{1}{c|}{\textbf{Case 1}} & \multicolumn{1}{c|}{1} & \multicolumn{1}{c|}{19.30} & \multicolumn{1}{c}{275.98} \\ 
    \multicolumn{1}{l|}{(correlation neglected)} & \multicolumn{1}{c|}{2} & \multicolumn{1}{c|}{24.61} & \multicolumn{1}{c}{273.01} \\ \hline
    \multicolumn{1}{c|}{\textbf{Case 2}}  & \multicolumn{1}{c|}{1} & \multicolumn{1}{c|}{0.99} & \multicolumn{1}{c}{343.07} \\ 
    \multicolumn{1}{l|}{(correlation considered)} & \multicolumn{1}{c|}{2} & \multicolumn{1}{c|}{1.99}  & \multicolumn{1}{c}{328.02} \\ \hline
    \end{tabular}
    \label{tab:results_correlation}
\end{table}

Neglecting the correlation in forecasts would always lead to the same performance as shown in Table~\ref{tab:results_correlation}: the level of risk desired would always be violated. This is due to the fact that correlated forecast errors tend to go in the same direction, therefore if any feeder has a high forecast error leading to very low net-load, it is likely that other feeders follow the same behaviour. Ignoring this compounded effect of forecast errors leads to an overly optimistic UFLS-allocation solution. This highlights the importance of computing the covariances of the forecast errors among feeders, and using these as inputs to the chance-constrained optimisation. The results in Table~\ref{tab:results_correlation} demonstrate that the proposed risk-constrained formulation for UFLS allocation can successfully accommodate correlated forecasts, meeting the desired level of overall risk.

\subsubsection{Advantages and limitations of assuming Gaussian uncertainty}\label{sec:results_GA}

In this section, the implications of solving the load shedding allocation using the Gaussian-assumption convex reformulation of the chance constraint are analysed. The optimisation problem in~(\ref{eq:stochastic_optimisation}) is solved by substituting constraint~(\ref{eq:cc}) by (\ref{gauss final form}), and the results are evaluated by comparing the performance of the selected feeders when the underlying uncertainty actually follows a different distribution. Gumbel, Laplace and Student's \textit{t} distributions are considered, while in all cases the first two moments are the ones defined in Table~\ref{tab:feedersdata}. 

\begin{figure}[H]
\vspace*{-5mm}
\hspace*{-8mm}
    \includegraphics[width=1.1\textwidth]{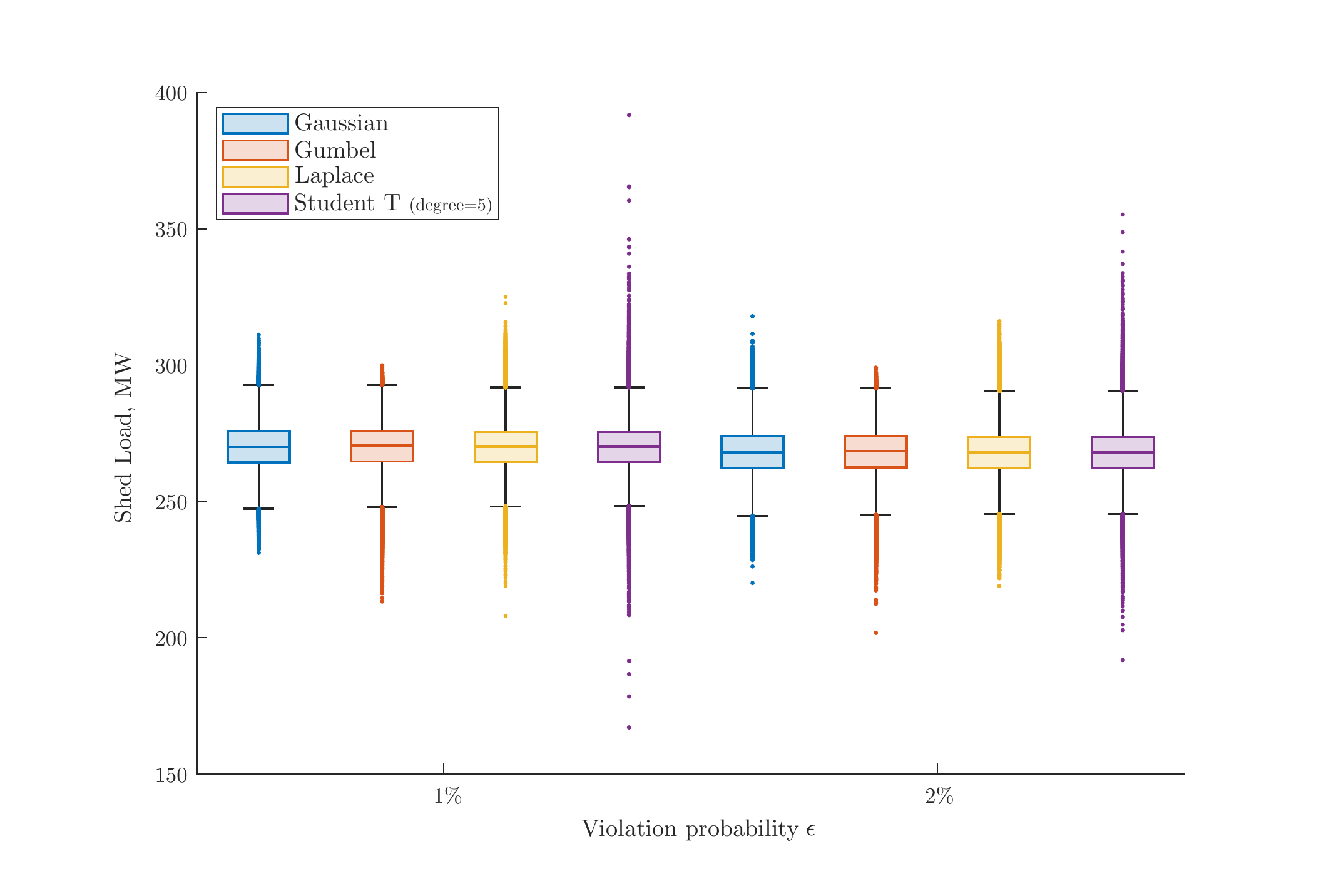}
    \vspace*{-5mm}
    \caption{Sampled distributions of available load in the selected feeders, for the four PDFs considered for underlying net-load uncertainty. Feeders were selected through a chance-constrained optimisation assuming Gaussian uncertainty.
    }
    \label{fig:12pCCGA}
\end{figure}

Table \ref{tab:results_ga} shows the results for two different risk thresholds, $\epsilon=1\%$ and $\epsilon=2\%$. The risk of shedding less load than required does indeed remain below $1\%$ or $2\%$, when uncertainty follows a Gaussian distribution. The advantage of constraint~(\ref{gauss final form}) is that it provides the least conservative feasible set for the UFLS allocation, if uncertainty can indeed be accurately modelled with Gaussian PDFs. However, when the underlying uncertainty in net-load through feeders follows a different PDF, the reliability requirements are not met. The highest deviation from the desired risk level is obtained when uncertainty follows a Gumbel distribution, which yields a $1.66\%$ risk for a desired value of under $1\%$. If the distribution is heavily tailed to both extremes, or to the left extreme, as shown in Figure \ref{fig:12pCCGA} for Laplace and Student's \textit{t}, and Gumbel distributions, respectively, the probability of not shedding the minimum required load would be higher than acceptable. These distributions imply that feeders are more likely to be consuming little power when called upon to participate in the UFLS scheme.

\begin{table}[h!]
\small
\caption{Results of chance-constrained allocation of load shedding assuming Gaussian uncertainty.}
\centering
    \begin{tabular}{c|l|c|c}
    \multicolumn{1}{l|}{$\epsilon$ (\%)} & \multicolumn{1}{c|}{\begin{tabular}[c]{@{}c@{}}Underlying \\ distribution\end{tabular}} & \multicolumn{1}{c|}{\begin{tabular}[c]{@{}c@{}}Samples below   \\ disconnection\\ requirement (\%)\end{tabular}} & \multicolumn{1}{c}{\begin{tabular}[c]{@{}c@{}}Expected load\\ disconnection (MW)\end{tabular}} \\ \hline
    \multirow{4}{*}{1} & Gaussian & 0.90 &  269.98 \\
     & Gumbel & 1.65 &  270.55 \\ 
     & Laplace & 1.11 &  270.01 \\ 
     & Student's \textit{t} & 1.08 &  270.01 \\ \hline
    \multirow{4}{*}{2} & Gaussian & 1.94 &  267.99 \\
     & Gumbel & 2.90 &  268.54 \\ 
     & Laplace & 2.13 &  267.99 \\ 
     & Student's \textit{t} & 2.07 &  268.00 \\ \hline
    \end{tabular}
    \label{tab:results_ga}
\end{table}

\subsubsection{Distributionally-robust chance-constrained optimisation}\label{sec:results_dr}

The limitations of the Gaussian chance constraint shown in the previous section can be overcome by using instead the distributionally robust chance constraint~(\ref{DRO final form}). The results of solving the UFLS-allocation problem with this constraint are shown in Table \ref{tab:results_dr}, which considers the same types of underlying uncertainty as in the previous section.

\begin{table}[h!]
\small
\caption{Results of distributionally-robust chance-constrained allocation of load shedding under different types of uncertainty.}
\centering
    \begin{tabular}{cccc}
    \multicolumn{1}{c|}{$\epsilon$ (\%)} & \multicolumn{1}{c|}{\begin{tabular}[c]{@{}c@{}}Underlying   \\ distribution\end{tabular}} & \multicolumn{1}{c|}{\begin{tabular}[c]{@{}c@{}}Samples below\\ disconnection \\ requirement (\%)\end{tabular}} &  \multicolumn{1}{c}{\begin{tabular}[c]{@{}c@{}}Expected load\\ disconnection (MW)\end{tabular}} \\ \hline
    \multicolumn{1}{c|}{\multirow{4}{*}{1}} & \multicolumn{1}{l|}{Gaussian} & \multicolumn{1}{c|}{0.00} & \multicolumn{1}{c}{357.98} \\  
    \multicolumn{1}{c|}{} & \multicolumn{1}{l|}{Gumbel} & \multicolumn{1}{c|}{0.00} &  \multicolumn{1}{c}{358.65} \\ 
    \multicolumn{1}{c|}{} & \multicolumn{1}{l|}{Laplace} & \multicolumn{1}{c|}{0.00} &  \multicolumn{1}{c}{357.98} \\ 
    \multicolumn{1}{c|}{} & \multicolumn{1}{l|}{Student's \textit{t}} & \multicolumn{1}{c|}{0.00} &  \multicolumn{1}{c}{358.01} \\ \hline
    \multicolumn{1}{c|}{\multirow{4}{*}{2}} & \multicolumn{1}{l|}{Gaussian} & \multicolumn{1}{c|}{0.00} & \multicolumn{1}{c}{317.99} \\ 
    \multicolumn{1}{c|}{} & \multicolumn{1}{l|}{Gumbel} & \multicolumn{1}{c|}{0.00} & \multicolumn{1}{c}{318.59} \\ 
    \multicolumn{1}{c|}{} & \multicolumn{1}{l|}{Laplace} & \multicolumn{1}{c|}{0.00} & \multicolumn{1}{c}{318.00} \\ 
    \multicolumn{1}{c|}{} & \multicolumn{1}{l|}{Student's \textit{t}} & \multicolumn{1}{c|}{0.00} & \multicolumn{1}{c}{318.00} \\ \hline
    \end{tabular}
\label{tab:results_dr}
\end{table}

\begin{figure}[h!]
\vspace*{-12mm}
\hspace*{-8mm}
    \centering
    \includegraphics[width=1.1\textwidth]{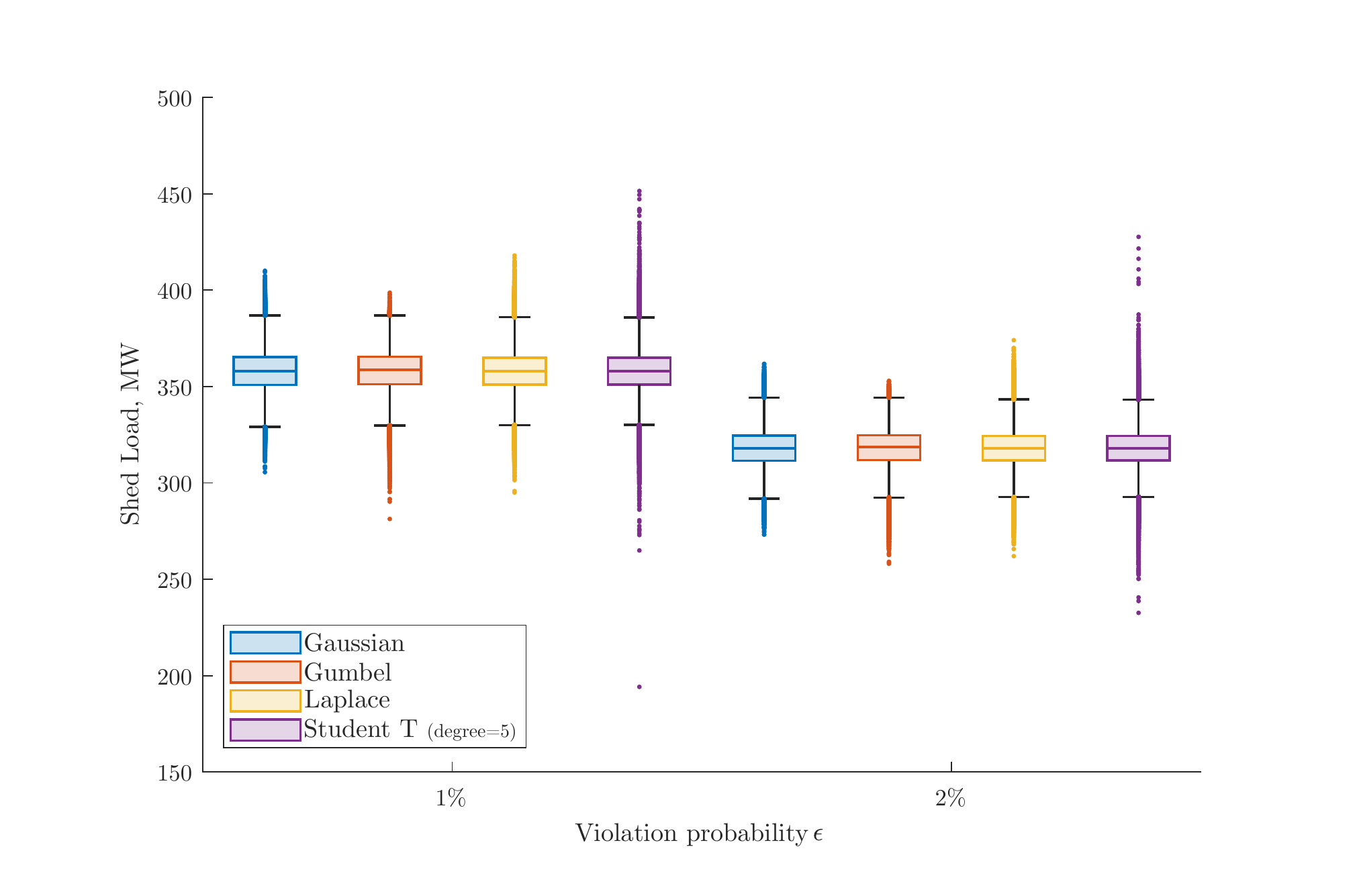}
    \vspace*{-12mm}
    \caption{Sampled distributions of available load in the selected feeders, for the four PDFs considered for underlying net-load uncertainty. Feeders were selected through a distributionally-robust chance-constrained optimisation.}
    \label{fig:12pCCDR}
\end{figure}

As expected, the constraint is robust against difference distributions and thus the reliability requirements are always met, i.e.,~the sampled sum of available load in all selected feeders is above the required load disconnection $L$ in most observations. In fact, the risk of violating the disconnection requirement is roughly zero for all these four PDFs considered, given that the distributionally-robust chance constraint is suitable for any PDF with these same values for the first two moments, and therefore the solution of the optimal UFLS allocation is conservative. Given that future distribution systems are expected to host heterogeneous load profiles, arising from a combination of DG, electric vehicles, heat pumps, combined heat and power plants, and other devices, the significant advantage of the distributionally-robust allocation lies on its ability to accommodate different types of uncertainty. Furthermore, the exact PDFs for net-load through each feeder are not needed, as simply the mean and standard deviation are required for applying the distributionally-robust method.

%%%%%%%%%%% FAIRNESS %%%%%%%%%%% 

\subsection{Introducing fairness measures in the UFLS-allocation process}

Finally, we discuss the issue of fairness in the selection of feeders that would potentially be disconnected through the UFLS scheme. Feeders with low penetration of DG would be more frequently chosen by the UFLS-allocation method, as they are less risky options due to their lower uncertainty. To address this issue and yield fairer results, we propose to include some degree of `synthetic uncertainty' when modelling these feeders in the risk-aware optimisation problem.

This synthetic uncertainty can be included through higher values of standard deviation. Before solving the optimisation problem to find the optimal selection of feeders, feeders that show significantly lower uncertainty than others must be identified, so that synthetic uncertainty can be introduced only in these feeders. While this modification would potentially result in more conservative load shedding, synthetic uncertainty will however not jeopardise compliance with reliability requirements imposed by the chance constraints. Note that the amount of synthetic uncertainty to be introduced is a design parameter, which must be tweaked to achieve a balance between fairness and deviating as little as possible from the truly optimal (yet socially biased) solution.

\begin{figure}[h!]
\vspace*{-5mm}
    \centering
    \includegraphics[width=0.9\textwidth]{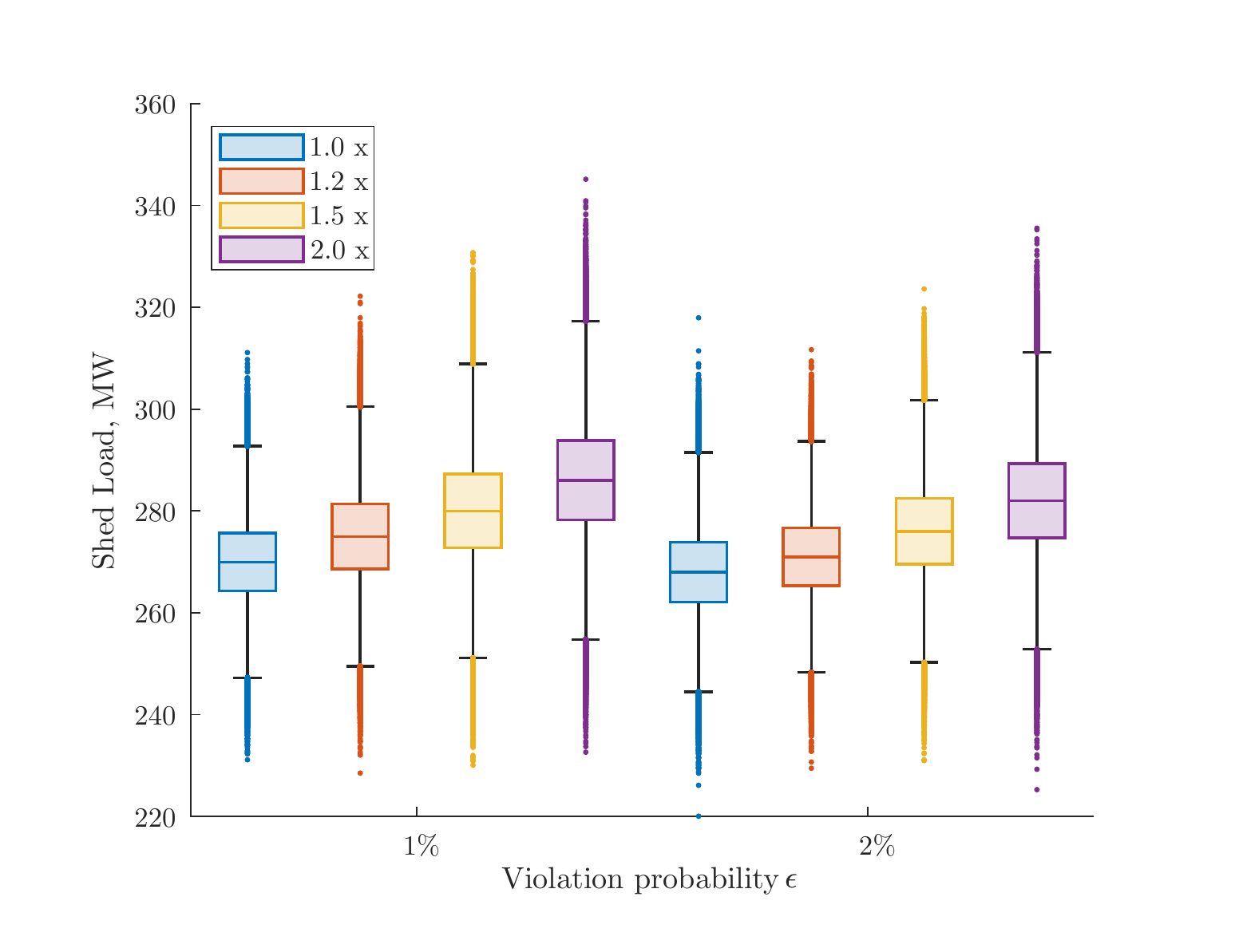}
    \vspace*{-5mm}
    \caption{Impact of introducing fairness measures in the UFLS-allocation process.}
    \label{fig:12pFAIRNESS}
\end{figure}

Figure~\ref{fig:12pFAIRNESS} shows the results of the allocation process for the base case optimisation (no synthetic uncertainty, tagged `1.0x') and the $20\%$ (1.2x), $50\%$ (1.5x) and $100\%$ (2.0x) increase in standard deviation for frequently selected feeders. Given that feeders chosen by the algorithm will change throughout time, depending on the exact forecast of net-load, the system operator could identify the feeders that are most frequently chosen to contribute towards the UFLS scheme, and include synthetic uncertainty only in these feeders. Here, the standard deviation of net-load for fourteen feeders, i.e.~all feeders except 1, 3, 8, 10, 14 and 15, was increased by the percentages detailed above. 

\begin{figure}[h!]
\vspace*{5mm}
    \centering
    \includegraphics[width=1\textwidth]{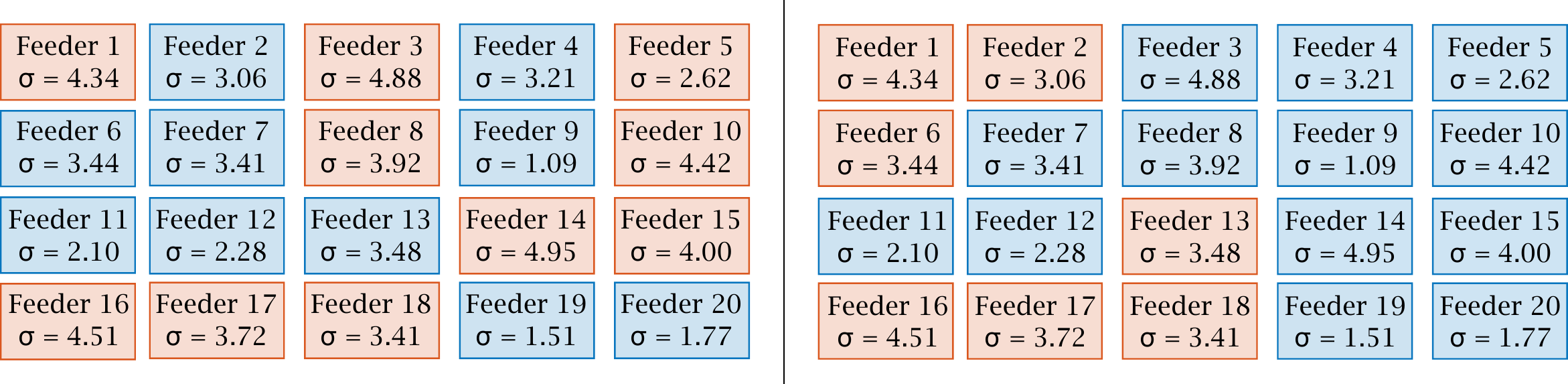}
    \caption{Change in selected feeders when introducing fairness measures in the UFLS allocation (orange = not selected, blue = selected). Left diagram shows the original allocation (without synthetic uncertainty), while right shows the socially-fair allocation.}
    \label{fig:rotation_fairness}
    \vspace*{5mm}
\end{figure}

The effect of synthetic uncertainty is to shift the selection of feeders contributing towards the UFLS scheme to achieve a fairer solution, as illustrated in Figure \ref{fig:rotation_fairness}. This figure shows an example of the change in the feeders' selection, when comparing the optimisation solutions of `1.0x' and `2.0x', for $\epsilon=1\%$. For the `1.0x' selection (right-hand side), feeders with low uncertainty (i.e.,~small standard deviation) are the ones chosen. Conversely, by introducing synthetic uncertainty, the selection of feeders on the right-hand side includes both truly uncertain feeders and feeders with low real uncertainty. Most importantly, the desired risk level is still met by the synthetic uncertainty introduced in the selection.

\section{Conclusion and future work} \label{sec:Conclusions}

In this paper, a risk-constrained methodology for selecting feeders contributing to the UFLS scheme has been proposed. Convex formulations for a Gaussian chance constraint and distributionally-robust chance constraint have been introduced for this purpose, enabling the system operator to explicitly consider uncertainty in net-load of each feeder during the UFLS-allocation process. This method is applicable both for systems with and without a communication network in place: for the latter, due to the inability to update UFLS relays close to real-time, considering uncertainty in net-load allows to find the optimal combination of
feeders to be chosen for the UFLS scheme, without risking under-delivery; while for the former (i.e.~WAMS-based adaptive UFLS schemes), the set-point for relays found through the proposed UFLS-allocation methodology can serve as a security layer to be used in case the communication network fails. 

Several relevant case studies have been conducted, demonstrating the importance of accounting for correlation in net-load of different feeders (notably due to correlated weather forecasts that affect the DG power output). Neglecting this positive correlation has been shown to be a risky practice, since the probability of falling short in total load disconnection would be higher than expected. 
The proposed chance-constrained formulation can also accommodate varying degrees of information on the characterisation of uncertainty, as a distributionally-robust framework has been deduced, which requires knowledge of only the two first moments of the underlying random variables.
Furthermore, a measure to increase fairness in the UFLS-allocation process has been introduced, since low-risk feeders (i.e.,~feeders with low standard deviation in net-load) would be more frequently chosen to contribute towards UFLS. Customers lacking DG would potentially be more frequently affected by interruption of energy supply, an issue that can be addressed by introducing some degree of `synthetic uncertainty' for these feeders.

As future work, the impacts of UFLS on other types of stability should be analysed, notably voltage stability. Since DG are already providing reactive power services to the transmission system in some countries, an unintended consequence of the activation of UFLS could be a lack/excess of reactive power. Therefore, the UFLS-allocation algorithm should be enhanced to include voltage stability metrics, avoiding an inadequate reactive power flow in the post-UFLS stage. Moreover, analysing the effect of varying cost for different feeders within the UFLS-allocation optimisation should be studied, given that some loads might be willing to participate in UFLS contracts in exchange of a fixed fee, while other loads would be subject to the value-of-lost-load penalty if eventually disconnected.

\vspace*{15mm}

\bibliographystyle{elsarticle-num.bst}
\bibliography{references.bib}
\end{document}